\documentclass[12pt,article]{amsart}
\usepackage{amsfonts}
\usepackage{amsthm}
\usepackage{amsmath}
\usepackage{amssymb} %
\usepackage{slashed}
\usepackage{amscd}
\usepackage[latin2]{inputenc}
\usepackage{t1enc}
\usepackage[mathscr]{eucal}
\usepackage{indentfirst}
\usepackage{graphicx}
\usepackage{graphics}
\usepackage{pict2e}
\usepackage{epic}
\numberwithin{equation}{section}
\usepackage[margin=3.5cm]{geometry}
\usepackage{epstopdf} 
\usepackage[utf8]{luainputenc}

\usepackage{tikz-cd}
\usepackage{dsfont}
\usetikzlibrary{matrix}

\newcommand\CP{{\mathbb C}\mathbb{P}}
\newcommand\QQ{{\mathbb Q}}
\newcommand\RR{{\mathbb R}}
\newcommand\ZZ{{\mathbb Z}}

\theoremstyle{plain}
\newtheorem{theorem}{Theorem}[section]
\newtheorem*{theorem*}{Theorem}

\newtheorem{corollary}[theorem]{Corollary}

 \theoremstyle{definition}

\newtheorem*{Conj*}{Conjecture}
\newtheorem{Conj}[theorem]{Conjecture}

\newtheorem*{problem*}{Problem}
\newtheorem{?}[theorem]{Problem}
\newtheorem{Question}[theorem]{Question}
\newtheorem*{Question*}{Question}

\newtheorem*{Acknowledgments}{Acknowledgments}

\begin{document}

\title[ ]{On the realization of symplectic algebras and rational homotopy types by closed symplectic manifolds}

\pagestyle{plain} 

\author{Aleksandar Milivojevi\'c}

\begin{abstract}  We answer a question of Oprea--Tralle on the realizability of symplectic algebras by symplectic manifolds in dimensions divisible by four, along with a question of Lupton--Oprea in all even dimensions. This will also allow us to address, in all even dimensions six and higher, another question of Oprea--Tralle on the possibility of algebraic conditions on the rational homotopy minimal model of a closed smooth manifold implying the existence of a symplectic structure on the manifold.
\end{abstract}

\address{Stony Brook University \\ Department of Mathematics} 

\email{aleksandar.milivojevic@stonybrook.edu}

\maketitle

\section{Introduction}

%{\let\thefootnote\relax\footnote{{2020 \textit{Mathematics Subject Classification}. Primary 55P62, 53C15.}}}
At the end of his famous two-page paper providing an example of a symplectic non-K\"ahler compact 4--manifold, Thurston \cite{Th76} posed the following conjecture:

\begin{Conj}\label{Conj1.1} (\cite{Th76}) Every closed $2k$--manifold which has an almost complex structure $\tau$ and a degree two real cohomology class $\alpha$ such that $\alpha^k \neq 0$ has a symplectic structure realizing $\tau$ and $\alpha$. \end{Conj}

Due to foundational results of Taubes and Witten in Seiberg--Witten theory, one can find counterexamples to this conjecture in dimension 4 (the argument to follow is well-known, see e.g. \cite[Example p.49]{Gom01}). Indeed, the oriented connected sum $\#_{i=1}^{2\ell+1} \CP^2$ for any $\ell \geq 1$ contains elements in $H^2$ not squaring to zero and admits an almost complex structure compatible with the orientation, but does not admit a compatible symplectic structure. By a classical result of Wu, one knows that a closed oriented four--manifold $M$ admits an almost complex structure if and only if there is a class $c \in H^2(M;\ZZ)$ such that its reduction mod 2 is the second Stiefel--Whitney class $w_2$ and $\int_M c^2 = 2\chi + 3\sigma$, where $\chi$ is the Euler characteristic and $\sigma$ is the signature. Since $w_2(\#_{i=1}^{2\ell+1}\CP^2) = (1,1,\ldots, 1) \in \ZZ_2^{\oplus 2\ell+1} \cong H^2(\#_{i=1}^{2\ell+1}\CP^2; \ZZ_2)$ and $2\chi + 3\sigma = 10\ell+9$, we see that $c = (3,1,3,1,\ldots, 1,3) \in H^2(\#_{i=1}^{2\ell+1}\CP^2;\ZZ)$ satisfies these conditions. Now, if $\#_{i=1}^{2\ell+1}\CP^2$ were to admit a symplectic structure realizing this almost complex structure, then by a theorem of Taubes, since $b_2^+ > 1$ (here $b_2^+$ is the dimension of the positive-definite subspace of the intersection form), it would have a non-vanishing Seiberg--Witten invariant. However, due to Witten, if a manifold is a connected sum of manifolds each with $b_2^+ \geq 1$ , then the Seiberg--Witten map is identically zero (see \cite[Corollary 4.1(2)]{Ko95}).

In dimensions $\geq 6$, these arguments from Seiberg--Witten theory do not directly apply, and Conjecture \ref{Conj1.1} remains open.

\vspace{0.7em}

In this note we will address the following variations of this conjecture:

\begin{Conj}\label{Conj1.2} (\cite[\textsection{6.5} Conjecture 3]{OT06}, \cite{HT08}, \cite{Tr00}) For every symplectic algebra $H$ over $\RR$, there is a closed symplectic manifold $M$ such that $H^*(M;\RR) \cong H$. \end{Conj}

A \textit{Poincar\'e duality algebra} (over the field $\Bbbk = \QQ$ or $\RR$) of dimension $n$ is a finite-dimensional graded-commutative algebra $H$ over $\Bbbk$ such that $H^n \cong \Bbbk$ and the pairing $H^{\ast} \otimes H^{n - \ast} \to \Bbbk$ given by $\alpha\otimes \beta \mapsto \mu(\alpha \beta)$ is non-degenerate for some (and hence any) choice of non-zero element $\mu \in {(H^n)}^{\ast}$. By a \textit{symplectic algebra} we mean a Poincar\'e duality algebra of dimension $2k$ for which there exists an element $\alpha \in H^2$ such that $\alpha^k \neq 0$. Hence, for simplicity, the adjectives "Poincar\'e duality" and "symplectic" will indicate properties of an algebra, not additional structure; $H^*(M;\RR)\cong H$ in the above conjecture will mean isomorphism of algebras. In dimensions $n = 4k$, a choice of orientation class $\mu$ lets one consider the signature of the induced pairing on $H^{2k}$. The pairing with respect to $a\mu$ will have the same signature for $a>0$, and the opposite signature for $a<0$; thus the signature of a $4k$--dimensional Poincar\'e duality algebra is well-defined up to sign. 

\begin{Question}\label{Q1.3} (\cite[Remark 2.11]{LO94}) Does a manifold that has rational cohomology algebra a symplectic algebra admit a symplectic structure? \end{Question}

In line with our previous definition, by a manifold we mean a connected orientable closed smooth manifold without a choice of orientation; hence admitting a symplectic (or almost complex) structure means possessing a symplectic form (or almost complex structure) inducing one of the two possible orientations on the manifold. Manifolds with symplectic rational cohomology algebras are also known as \textit{cohomologically symplectic} (or \textit{c--symplectic}) \cite{Tr00}.

\begin{Question}\label{Q1.4}(\cite[\textsection{6.5} Problem 4]{OT06}, \cite{Tr00}) Are there algebraic conditions on the minimal model $(\mathcal{M}_M, d)$ of a compact manifold $M$ implying the existence of a symplectic structure on $M$? \end{Question}

To answer Conjecture \ref{Conj1.2} in dimensions that are multiples of four, we will use a restriction on the topology of closed almost complex manifolds due to Hirzebruch. For Question \ref{Q1.3}, we will employ simply connected rational homology spheres not admitting spin${}^c$ structures in dimensions greater than five. This will immediately imply a negative answer to Question \ref{Q1.4}, in dimensions six and greater, when restricted to simply connected manifolds. In the non-simply connected (or more generally, non-nilpotent) case, one must first decide on what is meant by a minimal model in the sense of rational homotopy. However, we observe that any such notion which is invariant under weak homotopy equivalence of rationalizations in the sense of Bousfield--Kan cannot detect the existence of a symplectic form on a given manifold.

\section{Some symplectic algebras not realized by closed symplectic manifolds}

We provide counterexamples to Conjecture \ref{Conj1.2} in dimensions of the form $4k$. Consider for example \[ H = H^*\left( (S^2)^{2k} \# \mbox{\Large$\#$}_{i=1}^j (S^1\times S^{4k-1});\RR \right) \] for odd $j$. Taking $\alpha$ to be the sum of the images of generators of $H^2(S^2;\RR)$ under the inclusion $$H^2(S^2;\RR) \hookrightarrow H^2((S^2)^{2k};\RR) \hookrightarrow H,$$ we see that $\alpha^{2k} \neq 0$, and so $H$ is a symplectic algebra. Note that the signature $\sigma$ of the realizing oriented manifold $$(S^2)^{2k} \# \mbox{\Large$\#$}_{i=1}^j (S^1\times S^{4k-1})$$ is 0, and so the signature of any oriented manifold $M$ with $H^*(M;\RR) \cong H$ is 0, as the signature of a Poincar\'e duality algebra (with respect to any orientation class) is invariant up to sign under algebra isomorphisms of Poincar\'e duality algebras. On the other hand, the Euler characteristic satisfies $\chi = 2^{2k} - 2j \equiv 2 \bmod 4$ as $j$ is odd. By \cite[p.777]{Hir87}, a closed almost complex $4k$--manifold with the induced orientation satisfies the congruence $\chi \equiv (-1)^k \sigma \bmod 4$, so we conclude that $H$ cannot be realized by an almost complex manifold; in particular it cannot be realized by a symplectic manifold. We emphasize that this conclusion depends only on the algebra $H$, and so we have the following:

\begin{theorem} There are symplectic algebras $H$ over $\RR$ in every dimension $4k$, $k\geq 1$, such that there is no closed symplectic manifold $M$ with $H^*(M;\RR) \cong H$. \end{theorem}

Note that the above examples (i.e. the manifolds realizing them), by taking coefficients in $\QQ$ instead of $\RR$, provide an answer in the negative to Question \ref{Q1.3} in dimensions that are multiples of four. 

Alternatively, we can answer Question \ref{Q1.3} negatively in all even dimensions $\geq 6$ as follows: consider the Wu manifold $W = SU(3)/SO(3)$ of dimension 5; this is a simply connected rational homology sphere which does not admit a spin${}^c$ structure (see e.g. \cite[Example, p.50]{Fr00}). We consider the product $S^1 \times W$ and the result of performing surgery on the $S^1$ embedded in this product, i.e. removing an $S^1 \times \mathbb{D}^5$ from $S^1\times W$ and attaching a $\mathbb{D}^2 \times S^4$ along the common boundary $S^1 \times S^5$ by the identity map. This procedure is known as \textit{spinning} the manifold $W$ \cite{Suc90}. 

The result of spinning a simply connected rational homology sphere of dimension $n$ is a simply connected manifold \cite[Theorem IV.1.5]{Brow72}, and it is also a rational homology sphere, of dimension $n+1$ \cite[Lemma 2.1]{Suc90}. The result of spinning a manifold admits a spin${}^c$ structure if and only if the original manifold does \cite[Proposition 2.4]{AM19}. The latter two statements follow by a Mayer--Vietoris consideration, along with the characterization of manifolds admitting spin${}^c$ structures as those for which the second Stiefel--Whitney class $w_2$ is the mod 2 reduction of an integral class.

By iterating the spinning procedure applied to the non--spin${}^c$ Wu manifold $W$, we can thus produce a simply connected rational homology sphere $M^n$ of any dimension $n \geq 5$ not admitting a spin${}^c$ structure. For even $n$ we can then take the connected sum $M^n \# \CP^{n/2}$ of this rational homology sphere with $\CP^{n/2}$ to obtain a cohomologically symplectic but not symplectic manifold:

\begin{theorem} There are cohomologically symplectic manifolds in all dimensions $2k$, $k\geq 2$, that do not admit a symplectic structure. \end{theorem}

\begin{proof} The tangent bundle of the connected sum of oriented smooth $n$--manifolds $X \# Y$ is stably isomorphic to the bundle $\pi_X^* TX \oplus \pi_Y^*TY$, where $X\# Y \xrightarrow{\pi_X} X$ and $X\# Y \xrightarrow{\pi_Y} Y$ are the collapse maps. The manifolds are orientable, and therefore $w_2(X\#Y) = \pi_X^*w_2(X) + \pi_Y^*w_2(Y)$. Since the dimension is greater than two, we likewise have $H^2(X\#Y;R) \cong \pi_X^* H^2(X;R) \oplus \pi_Y^*H^2(Y;R)$ for any coefficient group $R$. 

A manifold admits a spin${}^c$ structure if and only if its second Stiefel--Whitney class $w_2$ is the mod 2 reduction of an integral class, and so from the above we see that $X\#Y$ is spin${}^c$ if and only if $X$ and $Y$ are spin${}^c$. (Note that the condition of admitting a spin${}^c$ structure does not depend on the choice of orientation.) An almost complex manifold (more generally, a stably almost complex manifold) admits a spin${}^c$ structure as the first Chern class of the stable almost complex structure reduces mod 2 to $w_2$. 

Now, in dimensions $n = 2k \geq 6$, taking the rational homology spheres $M^n$ previously constructed, which do not admit spin${}^c$ structures, we see that $M^n \# \CP^{n/2}$ is a cohomologically symplectic manifold which does not admit a spin${}^c$ structure, and hence does not admit an almost complex structure; in particular it does not admit a symplectic structure. Together with the discussion in the introduction, this covers all dimensions $\geq 4$. \end{proof}

As far as Question \ref{Q1.3} is concerned, we could have replaced $M^n$ above with any non--spin${}^c$ orientable closed manifold (we remark that the smallest dimension in which these exist is 5). However, having non--spin${}^c$ simply connected rational homology spheres in hand will let us immediately turn to Question \ref{Q1.4} in the following section.

\section{The existence of a symplectic structure cannot be detected from the rational homotopy model}

We now address Question \ref{Q1.4}. For any simply connected symplectic manifold $X$ of dimension at least six, consider the connected sum $M \# X$ (to form the connected sum we choose any orientation on $M$ and $X$), where $M$ is a non-spin${}^c$ simply connected rational homology sphere as in the previous section. The collapse map $M \# X \to X$ is a rational homotopy equivalence as it is a rational homology equivalence of simply connected spaces, and so the minimal models of these manifolds are isomorphic while only one of them admits a symplectic structure (with respect to \textit{some} orientation), as $M\#X$ does not admit a spin${}^c$ structure. Since $X$ was an arbitrary simply connected symplectic manifold, we conclude that there can be no algebraic condition on minimal models of simply connected manifolds which implies the manifold admits a symplectic structure. 

In the non-simply connected case, the classical theory for simply connected spaces of finite type due to Sullivan extends immediately to spaces with nilpotent fundamental group which acts nilpotently on the higher homotopy groups, and the algebraic information encoded in the minimal model directly corresponds to geometric information. Bousfield and Kan extended the procedure of rationalizing spaces to all path-connected spaces in two ways \cite{BoKa71}: the $\QQ$--completion and the fiberwise $\QQ$--completion, both restricting to the classical rationalization on nilpotent spaces (see \cite{RWZ19} for an overview). A map $X \to Y$ induces a weak homotopy equivalence of $\QQ$--completions if it induces an isomorphism on rational homology \cite{BoKa71}, and it induces a weak homotopy equivalence of fiberwise $\QQ$--completions if it induces an isomorphism on fundamental groups and on rationalized higher homotopy groups (see \cite[Theorem 3]{RWZ19}). Substantial progress has been made in algebraically encoding spaces up to these notions of equivalence, extending the classical theory of rational homotopy minimal models; see \cite{GHT00}, \cite{BFMT18}. 

We now observe that for any (not necessarily simply connected) symplectic manifold $X$ of dimension $\geq 6$, there is another manifold, not admitting a symplectic structure, which is equivalent to $X$ under either of the above notions. Consider again the collapse map $M\# X \to X$, where $M$ is a non-spin${}^c$ simply connected rational homology sphere; this map induces an isomorphism on rational homology and hence a weak homotopy equivalence of $\QQ$--completions. 

The collapse map likewise induces an isomorphism of fundamental groups, and to verify it induces an isomorphism on $\pi_{\geq 2} \otimes \QQ$, we proceed as follows: pick basepoints and consider the induced map $$\widetilde{M \# X} \to \widetilde{X}$$ on universal covers. Since $M$ is simply connected, the space $\widetilde{M \# X}$ can be visualized as the universal cover of $X$, with a small disk $\mathbb{D}_i$ around each preimage $\tilde{b}_i$ of the basepoint of $X$ (chosen to coincide with the center of the disk at which the connected sum with $M$ is performed) replaced by $M \# \mathbb{D}_i$. The map on universal covers $\widetilde{M \# X} \to \widetilde{X}$ is then the collapse map $M \# \mathbb{D}_i \to \mathbb{D}_i$ applied at each of these disks, and the identity elsewhere; see Figure 1.

Consider the open cover of $\widetilde{M\# X}$ given by a small neighborhood of $\bigcup_i (M \# \mathbb{D}_i)$ and the complement of $\bigcup_i (M \# \mathbb{D}_i)$, along with the open cover of $\widetilde{X}$ given by a small neighborhood of $\bigcup_i \mathbb{D}_i$ and the complement of $\bigcup_i \mathbb{D}_i$. 

Applying the naturality of the Mayer--Vietoris sequence in homology to these open covers, by the five lemma we see that the map $\widetilde{M \# X} \to \widetilde{X}$ induces an isomorphism on rational homology. Indeed, $M \# \mathbb{D}_i$ has the homotopy type of $M$ with a point removed, and so has the rational homology groups of a point; thus the map $\bigcup_i (M \# \mathbb{D}_i) \to \bigcup_i \mathbb{D}_i$ induces an isomorphism on rational homology. On the complement and the intersection the map is the identity.

\tikzset{every picture/.style={line width=0.75pt}}

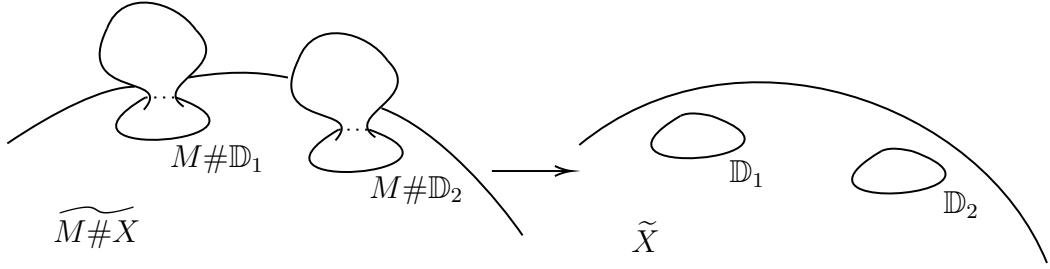
\begin{figure}\label{map}
\begin{tikzpicture}[x=0.58pt,y=0.58pt,yscale=-1,xscale=1]
%uncomment if require: \path (0,377); %set diagram left start at 0, and has height of 377

%Curve Lines [id:da9503134268347189] 
\draw    (110,122) .. controls (131,98) and (63,116) .. (86,69) ;
%Curve Lines [id:da7635299332709866] 
\draw    (136,120) .. controls (114,105) and (169,95) .. (144,70) ;
%Curve Lines [id:da7554818896928854] 
\draw    (86,69) .. controls (96,42) and (131,51) .. (144,70) ;
%Straight Lines [id:da031221125583767773] 
\draw    (337,162) -- (355,162) -- (387,162) ;
\draw [shift={(389,162)}, rotate = 180] [color={rgb, 255:red, 0; green, 0; blue, 0 }  ][line width=0.75]    (10.93,-3.29) .. controls (6.95,-1.4) and (3.31,-0.3) .. (0,0) .. controls (3.31,0.3) and (6.95,1.4) .. (10.93,3.29)   ;
%Curve Lines [id:da08919677069664722] 
\draw    (111,114) .. controls (36,158) and (212,143) .. (131,114) ;
%Straight Lines [id:da3188999349318469] 
\draw  [dash pattern={on 0.84pt off 2.51pt}]  (111,114) -- (131,114) ;
%Curve Lines [id:da05591850258443243] 
\draw    (261,140) .. controls (239,125) and (294,115) .. (269,90) ;
%Curve Lines [id:da9455950670550153] 
\draw    (211,89) .. controls (221,62) and (256,71) .. (269,90) ;
%Curve Lines [id:da32932538050602855] 
\draw    (235,142) .. controls (256,118) and (188,136) .. (211,89) ;
%Curve Lines [id:da25403192301413635] 
\draw    (237,135) .. controls (162,179) and (338,164) .. (257,135) ;
%Straight Lines [id:da431765866957712] 
\draw  [dash pattern={on 0.84pt off 2.51pt}]  (237,135) -- (257,135) ;
%Curve Lines [id:da1657351893037189] 
\draw    (460,126) .. controls (385,170) and (561,155) .. (480,126) ;
%Curve Lines [id:da4529796563105093] 
\draw    (591,149) .. controls (538.99,179.51) and (607.69,181.65) .. (627.68,170.42) .. controls (636.51,165.46) and (635.83,157.89) .. (611,149) ;
%Curve Lines [id:da1841477457541687] 
\draw    (21,144) .. controls (52,123) and (84,107) .. (104,107) ;
%Curve Lines [id:da1758909975017171] 
\draw    (139,101) .. controls (174,95) and (193,99) .. (204,101) ;
%Curve Lines [id:da4241590408207234] 
\draw    (265,121) .. controls (303,135) and (342,182) .. (357,204) ;
%Curve Lines [id:da06156692447344203] 
\draw    (394,145) .. controls (496,61) and (656,116) .. (699,222) ;
%Curve Lines [id:da7557405619219453] 
\draw    (460,126) .. controls (466,123) and (475,125) .. (480,126) ;
%Curve Lines [id:da43645663158819725] 
\draw    (591,149) .. controls (597,146) and (606,148) .. (611,149) ;

% Text Node
%\draw (133,38) node [anchor=north west][inner sep=0.75pt]   [align=left] {$\displaystyle M$};
% Text Node
%\draw (259,57) node [anchor=north west][inner sep=0.75pt]   [align=left] {$\displaystyle M$};
% Text Node
\draw (49,183) node [anchor=north west][inner sep=0.75pt]   [align=left] {$\displaystyle \widetilde{M\#X}$};
% Text Node
\draw (426,193) node [anchor=north west][inner sep=0.75pt]   [align=left] {$\displaystyle \widetilde{X}$};
% Text Node
\draw (125,143) node [anchor=north west][inner sep=0.75pt]   [align=left] {$\displaystyle M \#\mathbb{D}_1$};
% Text Node
\draw (255,163) node [anchor=north west][inner sep=0.75pt]   [align=left] {$\displaystyle M\# \mathbb{D}_2$};
% Text Node
\draw (490,153) node [anchor=north west][inner sep=0.75pt]   [align=left] {$\displaystyle \mathbb{D}_1$};
% Text Node
\draw (629.68,173.42) node [anchor=north west][inner sep=0.75pt]   [align=left] {$\displaystyle \mathbb{D}_2$};

\end{tikzpicture}
\caption{The map between universal covers.}
\end{figure}

Thus, since $\widetilde{M \# X}$ and $\widetilde{X}$ are simply connected, this map between universal covers induces an isomorphism on rational homotopy groups. Now, from the map of long exact sequences of homotopy groups associated to the map of fibrations

$$\begin{tikzcd}
\widetilde{M \# X} \arrow[d] \arrow[r] & \widetilde{X} \arrow[d] \\
M\#X \arrow[r]                         & X                      
\end{tikzcd}$$

\noindent and the five lemma again, we conclude that $M \# X \to X$ induces an isomorphism on $\pi_{\geq 2} \otimes \QQ$. Therefore the fiberwise $\QQ$--completions of these spaces are also equivalent. In conclusion, we have:
\begin{theorem} There are no algebraic conditions on the minimal model $(\mathcal{M}_X, d)$ of a manifold $X$ implying the existence of a symplectic structure on $X$, in dimensions six or greater.\end{theorem}

Here by a minimal model we mean any object (in particular, the classical minimal models in the case of finite-type nilpotent spaces) which is invariant up to isomorphism under weak homotopy equivalence of rationalizations in either sense of Bousfield--Kan.

We note that the same argument, using non-spin${}^c$ simply connected rational homology spheres, shows that the existence of a (stable almost) complex structure cannot be implied by algebraic conditions on the minimal model:

\begin{corollary}  There are no algebraic conditions on the minimal model of a manifold $X$ implying the existence of a complex structure (or more generally a stable almost complex structure) on $X$, in dimensions six or greater. \end{corollary}

\section{Another variation of Thurston's conjecture}

It seems that the following question, another variation of Conjecture \ref{Conj1.1}, is still unanswered in all dimensions $\geq 4$:

\begin{Question} Is there a symplectic algebra which is realized by a closed almost complex manifold but not realized by a closed symplectic manifold? \end{Question}

Currently there are no known topological obstructions to a closed smooth manifold admitting a symplectic structure beyond those of admitting an almost complex structure and having a symplectic cohomology algebra. A possible direction presents itself as it seems that for all known examples of closed symplectic $2n$--manifolds, the Betti numbers $b_i$ for $i\leq n$ satisfy the non-decreasing property $b_0 \leq b_2 \leq b_4 \leq \cdots$ and $b_1 \leq b_3 \leq \cdots$ \cite[Question 1.1]{Cho16}. A proof that this property holds for all closed symplectic manifolds would immediately enable one to provide counterexamples to Conjectures \ref{Conj1.1} and \ref{Conj1.2}, along with the above question.

\begin{Acknowledgments} The author would like to thank Scott Wilson for his numerous helpful suggestions, and John Morgan for a relevant discussion on a separate problem.\end{Acknowledgments}

\end{document}